\documentclass[12pt,reqno,a4paper]{amsart}
\usepackage{amssymb,amscd}

\usepackage{hyperref}

\begin{document}

\let\kappa=\varkappa
\let\eps=\varepsilon
\let\phi=\varphi
\let\p\partial

\def\Z{\mathbb Z}
\def\R{\mathbb R}
\def\C{\mathbb C}
\def\Q{\mathbb Q}

\def\OO{\mathcal O}
\def\CP{\C{\mathrm P}}
\def\RP{\R{\mathrm P}}
\def\conj{\overline}
\def\Beta{\mathrm{B}}
\def\const{\mathrm{const}}

\def\H{\mathrm H}
\def\Sph{S}
\def\B{B}

\renewcommand{\Im}{\mathop{\mathrm{Im}}\nolimits}
\renewcommand{\Re}{\mathop{\mathrm{Re}}\nolimits}
\newcommand{\codim}{\mathop{\mathrm{codim}}\nolimits}
\newcommand{\id}{\mathop{\mathrm{id}}\nolimits}
\newcommand{\Aut}{\mathop{\mathrm{Aut}}\nolimits}
\newcommand{\lk}{\mathop{\mathrm{lk}}\nolimits}
\newcommand{\Ker}{\mathop{\mathrm{Ker}}\nolimits}
\newcommand{\sign}{\mathop{\mathrm{sign}}\nolimits}
\newcommand{\rk}{\mathop{\mathrm{rk}}\nolimits}

\renewcommand{\mod}{\mathrel{\mathrm{mod}}}

\newtheorem*{mainthm}{Theorem}
\newtheorem{thm}{Theorem}
\newtheorem{lem}[thm]{Lemma}
\newtheorem{prop}[thm]{Proposition}
\newtheorem{cor}[thm]{Corollary}

\theoremstyle{definition}
\newtheorem{exm}[thm]{Example}
\newtheorem{rem}[thm]{Remark}
\newtheorem{df}[thm]{Definition}

\renewcommand{\thesubsection}{\arabic{subsection}}

\title{Lagrangian Klein bottles in $\R^{2n}$}
\author{Stefan Nemirovski}
\thanks{The author was supported by grants from DFG, RFBR, Russian Science Support Foundation,
and the programme ``Leading Scientific Schools of Russia.''}
\address{%
Steklov Mathematical Institute;\hfill\break
\strut\hspace{8 true pt} Ruhr-Universit\"at Bochum}
\email{stefan@mi.ras.ru}


\maketitle

The $n$-dimensional Klein bottle $K^n$, $n\ge 2$, is obtained by gluing the
ends of the cylinder $\Sph^{n-1}\times [0,1]$ via an orientation reversing
isometry of the standard $(n-1)$-sphere $\Sph^{n-1}\subset\R^n$.

\begin{mainthm}
\label{main}
The $n$-dimensional Klein bottle $K^n$ admits a Lagrangian embedding into
the standard symplectic $2n$-space $(\R^{2n},\omega_0)$ if and only if $n$ is odd.
\end{mainthm}

The existence of Lagrangian embeddings of {\it odd\/}-dimensional Klein
bottles into $(\R^{2n},\omega_0)$ was proved by Lalonde~\cite{Lal}.
As observed in~\cite{N}, an explicit embedding is suggested by Picard--Lefschetz theory.
Indeed, the antipodal map
$$
\R^{2k+1}\supset\Sph^{2k}\ni (x_1,x_2,\dots,x_{2k+1})\longmapsto
(-x_1,-x_2,\dots,-x_{2k+1})\in \Sph^{2k}\subset\R^{2k+1}
$$
reverses the orientation on~$\Sph^{2k}\subset\R^{2k+1}$
and therefore the formula
\begin{equation}
\label{odd}
\Sph^{2k}\times[0,1]
\ni(x_1,\dots,x_{2k+1},t)\longmapsto
(e^{\pi i t}x_1,\dots,e^{\pi i t}x_{2k+1})\in\C^{2k+1}
\end{equation}
defines an embedding of the odd-dimensional Klein bottle $K^{2k+1}$ into~$\C^{2k+1}=\R^{4k+2}$.
It is easy to check that this embedding is Lagrangian with respect
to the standard symplectic form $\omega_0=\frac{i}{2}\sum dz_\ell\wedge d\conj z_\ell$.

Thus, the main task is to prove that an {\it even\/}-dimensional Klein bottle does not
admit a Lagrangian embedding into $(\R^{2n},\omega_0)$.
For the usual Klein bottle~$K^2$, this problem was proposed by Givental'~\cite{Gi} and
resolved by Shevchishin~\cite{Sh}. The argument in the general case follows the lines
of the author's alternative proof of Shevchishin's result~\cite{N1} (cf.\ also~\cite{EP}).
Namely, self-linking invariants introduced by Rokhlin and Viro
are used to show that a suitable Luttinger-type surgery along a Lagrangian
$K^{2k}\subset\R^{4k}$ would produce an impossible symplectic manifold.

\subsection{Rokhlin and Viro indices for totally real Klein bottles}
\label{ind}
Let us fix $n$ and denote the $n$-dimensional Klein bottle simply by $K$.
Let $m\subset K$ be a fibre of the natural fibre bundle $K\to \Sph^1$. Then $m$ is
an embedded $(n-1)$-dimensional sphere in $K$. Note that $m$ is co-orientable
and choose a non-vanishing normal vector field $\nu_{m,K}$ on~$m$.

Consider now a totally real embedding $K\hookrightarrow\C^n$, i.\,e., an embedding
such that $T_pK$ is transversal to $iT_pK$ at every point $p\in K$.
(For a Lagrangian embedding, the subspaces $T_pK$ and $iT_pK$ would be orthogonal
with respect to the standard metric on $\C^n$.)
Let $m^\sharp$ be the pushoff of $m$ in the direction of the vector field $i\nu_{m,K}$.
The $\mod 2$ homology class $[m^\sharp]\in \H_{n-1}(\C^n\setminus K;\Z/2)$ is independent of the
choice of~$\nu_{m,K}$. The linking number
$$
V=\lk(K,m^\sharp)\in\Z/2
$$
is called the {\it Viro index\/} of $m\subset K$ (cf.\ \cite{N1}, \S 1.2).

In order to compute $V$, we choose an immersed $n$-ball $M=\iota(\B^n)\subset\C^n$ such that
\begin{itemize}
\item[a)] $\p M=\iota(\p\B^n)=m$, and $M$ is normal to $K$ along $m$;
\item[b)] the self-intersections of $M$ and the intersections of its interior with $K$ are transverse double points;
\item[c)] the tangent (half-)space of $M$ at a point $p\in m=\p M$ is spanned
by $T_p m$ and $i\nu_{m,K}$;
\item[d)] the $\R\C$-singular points of $M$ are generic (see the definitions in \S\ref{rc} below).
\end{itemize}
(Immersions satisfying (a) and (b) are called {\it membranes\/} spanned by~$m$.)
Note that by (c) the pushoff of $m$ inside $M$ is precisely $m^\sharp$, and hence
\begin{equation}
\label{lk}
V=\#(M\cap K)\mod 2,
\end{equation}
where $\#(M\cap K)$ denotes the number of interior intersection points of $M$ and $K$.

Suppose now that $n$ is even and consider the {\it Rokhlin index\/} of $M$ defined by the formula
\begin{equation}
\label{rdef}
R=n(M,\nu_{m,K}) + \# (M\cap K),
\end{equation}
where $n(M,\nu_{m,K})\in \Z$ is the obstruction to extending $\nu_{m,K}$ to a non-vanishing normal vector field on~$M$,
i.\,e., the algebraic number of zeroes of a generic normal extension. (For odd $n$, this number
is defined only $\mod 2$.)

\begin{lem}
\label{rokh}
$R=0\mod 2$.
\end{lem}

This will be proved in \S\ref{rproof} using nothing much. Note, however, that this
is the only place where the assumption that $n$ is even will be used in a crucial way
(see Remark~\ref{caution}).

\begin{lem}
\label{norm}
$n(M,\nu_{m,K})=1\mod 2$.
\end{lem}

This will be proved in \S\ref{nproof} using a topological count of $\R\C$-singularities
recalled briefly in \S\ref{rc} following Domrin~\cite{Do}.

\begin{lem}
\label{viro}
$V=1\mod 2$.
\end{lem}

This follows immediately from  formulas~(\ref{lk}) and~(\ref{rdef}) and the preceding two lemmas and
will play a key role in the proof of the main theorem in \S\ref{surgery} and~\S\ref{mainproof}.

\subsection{Proof of Lemma~\ref{rokh}}
\label{rproof}
Cut $K$ along $m$ and glue two copies of $M$ into the resulting `holes'
to obtain an $n$-sphere $S$.  Choose an orientation on $S$ and note that
it induces the {\it same\/} orientation on each of the two copies of $M$.
(If we had $\Sph^{n-1}\times\Sph^1$ instead of the Klein bottle, the orientations
would be opposite.) Let $\nu$ be a generic normal extension of $\nu_{m,K}$ to~$M$.
Transform $S$ into a generically immersed sphere by pushing the two copies of $M$
apart in the direction of $\nu$ and then smoothing the result.

Now we can compute the normal Euler number of $S$ and the algebraic number
of its double points. Namely,
$$
n(S)=n(K)+2n(M,\nu_{m,K})=2n(M,\nu_{m,K}),
$$
where we have used the fact that for a totally real embedding of $K$
the normal Euler number is equal to the Euler characteristic of $K$
which is zero. Similarly,
$$
\#_{alg}(S)=n(M,\nu_{m,K})+2\,\#_{alg}(M\cap K)+4\,\#_{alg}(M),
$$
where the signs in $\#_{alg}(M\cap K)$ and $\#_{alg}(M)$ are given by the induced orientations
on $M$ and $K$ as subsets of~$S$.

On the other hand, by the usual formula for the homological
self-intersection index of an {\it oriented\/} immersed submanifold,
we have
$$
[S]\cdot[S]=n(S)+2\#_{alg}(S)=4\bigl(n(M,\nu_{m,K})+\#_{alg}(M\cap K) +2\,\#_{alg}(M)\bigr).
$$
The homology class $[S]$ is obviously trivial in $\C^n$, hence
$$
n(M,\nu_{m,K})+\#_{alg}(M\cap K)+ 2\,\#_{alg}(M) =0,
$$
and the result follows from~(\ref{rdef}) because $\#_{alg}(M\cap K)=\#(M\cap K)\mod 2$. \qed

\begin{rem}
The above argument and the result for $n=2$ go back to Rokhlin
(see~\cite{GM} and the proof of Lemma~1.12 in~\cite{N1}). Note that we are
actually proving a congruence  modulo~8 using a trivial case of van der Blij's lemma to conclude
that $[S]\cdot[S]=0\mod 8$.
\end{rem}

\begin{rem}
\label{caution}
For odd $n$, the residue $R\mod 2$ is well-defined but the lemma is false. (Our proof does not work
because the intersection index is not symmetric.) For instance,
for the embedding given by~(\ref{odd}),
the totally real $n$-ball $\{(z_1,\ldots,z_n)\in\C^n\mid z_j\in\R, \|z\|\le 1\}$
is a membrane satisfying conditions (a)-(d) and such that $R=1\mod 2$.

\end{rem}

\subsection{$\R\C$-singularities and characteristic classes}
\label{rc}
Here's a digression needed for the proof of Lemma~\ref{norm}.
The material is mostly taken from~\cite{Do} (cf.\ also~\cite{La} and~\cite{We}).

Let $j:N\to\C^n$ be an immersion of a real oriented $n$-dimensional manifold. (Note that the real
dimension of $N$ is equal to the complex dimension of~$\C^n$.) A point $p\in N$ is called {\it $\R\C$-singular\/}
if the dimension of the maximal complex subspace in $j_*T_p N\subset\C^n$ is positive (i.\,e., larger than expected).
This dimension is called the {\it order\/} of an $\R\C$-singular point.
Denote by $C_\mu(N)$ the set of $\R\C$-singular points of order~$\mu$ and by $C(N)$ the set of all $\R\C$-singular points.

Let $j_*^\C: TN\otimes\C\to j^*T\C^n$ be the complex vector bundle map given by $j_*^\C(v\otimes \lambda):=\lambda j_*(v)$.
Its kernel at a point $p\in N$ is isomorphic to the maximal complex subspace of~$j_*T_pN$.
Thus, $C_\mu(N)$ coincides with the singularity set $\Sigma_\mu=\{p\in N\mid \rk_\C j_*^\C= n-\mu\}$.
If the immersion $j$ is generic, then by~\cite{Do}, Lemma~1.3, the bundle map $j_*^\C$
is generic in the sense of~\cite{McPh}.
Hence, each $C_\mu(N)$ is an oriented $(n-2\mu^2)$-dimensional submanifold, $C(N)=\conj{C_1(N)}$,
and there exists a canonical desingularisation $\widetilde\Sigma_1\to \conj\Sigma_1=C(N)$
such that the complex line bundle $\Ker j_*^\C|_{\Sigma_1}$ extends to $\widetilde\Sigma_1$.
(Explicitly, $\widetilde\Sigma_1$ is the closure of the image of $\Ker j_*^\C|_{\Sigma_1}$
in the projectivisation  of $TN\otimes\C$ and the extension of $\Ker j_*^\C|_{\Sigma_1}$
is given by the tautological line bundle.) The extended bundle lies in the kernel of
the pull-back of ${j_*^\C}$ to~$\widetilde\Sigma_1$ and therefore corresponds to
a complex line subbundle $\mathcal Z$ of the pull-back of~$TN$.

Assume now that the manifold $N$ is compact without boundary and its dimension $n$ is even.
Define $[C(N)]$ as the fundamental class of the oriented manifold~$\widetilde\Sigma_1$.
Then from Theorem~3 and Remark~1.4 in~\cite{Do} one obtains the formula
\begin{equation}
\label{domrin}
\langle c_1({\mathcal Z})^{(n/2-1)},[C(N)]\rangle = \langle c_{n/2}(-TN\otimes\C), [N]\rangle,
\end{equation}
where $-TN\otimes\C$ denotes the $K$-theoretic inverse of $TN\otimes\C$. In the statements of
the results in~\cite{Do} it is assumed that $C(N)=C_1(N)$ but the proofs carry over to the general
case with formal changes. ($\widetilde\Sigma_1$ has to be used instead of the set~$\Sigma$
introduced on p.\,910~of~\cite{Do}.)

\begin{rem}
\label{lai}
The $\R\C$-singular points of a generic immersed surface in $\C^2$ are isolated
(and more often referred to as `complex points' or `complex tangencies'). Formula (\ref{domrin})
reduces in this case to the elementary formula $I_+ - I_-=0$, where the {\it Lai indices\/}~$I_\pm$
of an oriented immersed surface are defined by counting its complex points  with suitable signs
(see~\cite{Do},~\S 3).
\end{rem}

\subsection{Proof of Lemma~\ref{norm} {\rm (cf.\ \cite{N1}, Proof of Lemma~1.13)}}
\label{nproof}
Let us construct a normal extension of $\nu_{m,K}$ to $M$ in the following way.
Consider the vector field $i\nu_{m,K}$. It is tangent to $M$ and transverse to $\p M$
by the choice of~$M$. Let $\tau$ be an extension of this vector field to a tangent
vector field on $M$ with a single transverse zero. (Recall that $M$ is a ball.)
Then $-i\tau$ gives a normal extension of $\nu_{m,K}$ that vanishes at the
zero of $\tau$ and at the points where $\tau$ lies in a non-trivial complex subspace
contained in $T_pM$. For a sufficiently generic~$\tau$, the latter points
lie in $C_1(M)$. In other words, we have to count the zeroes of a (generic) section
of the quotient bundle $\mathcal E=\widetilde{TM}/\mathcal Z$, where $\widetilde{TM}$
is the pull-back of $TM$ to $\widetilde\Sigma_1$.
As we only need the answer mod 2, it is given by the evaluation of the top Stiefel--Whitney
class $w_{n-2}(\mathcal E)$ on the fundamental class~$[C(M)]:=[\widetilde\Sigma_1]$.

Since $TM$ is trivial, we have
\begin{equation}
\label{wh}
1=(1+w_2(\mathcal Z))(1+w_1(\mathcal E)+\dots+w_{n-2}(\mathcal E))
\end{equation}
by the Whitney formula. It follows immediately that
$$
w_{n-2}(\mathcal E)=w_2(\mathcal Z)^{(n/2-1)}
$$
and hence
$$
\langle w_{n-2}(\mathcal E),[C(M)]\rangle = \langle w_2(\mathcal Z)^{(n/2-1)},[C(M)]\rangle =
\langle c_1({\mathcal Z})^{(n/2-1)},[C(M)]\rangle \mod 2.
$$
In order to show that the latter quantity vanishes (already as an integer), we apply formula~(\ref{domrin})
to an immersed sphere $S$ similar to the one used in the proof of Lemma~\ref{rokh} above.
Namely, we glue two copies of $M$ to $K$ cut along $m$ but this time only smoothen
the result near $m$. Condition (c) in \S\ref{ind} ensures that this smoothing can be done
so that no additional $\R\C$-singularities are created and hence the set $C(S)$ consists
of two copies of~$C(M)$ with the same orientation and the same line bundle $\mathcal Z$.
Thus,
$$
2\langle c_1({\mathcal Z})^{(n/2-1)},[C(M)]\rangle=\langle c_1({\mathcal Z})^{(n/2-1)},[C(S)]\rangle
\overset{\text{(\ref{domrin})}}{=}\langle c_{n/2}(-TS\otimes\C), [S]\rangle=0.
$$
It follows that $\langle w_{n-2}(\mathcal E),[C(M)]\rangle=0\mod 2$ and hence the normal projection
of $-i\tau$ has an odd number of zeroes, which proves that $n(M,\nu_{m,K})=1\mod 2$.\qed

\begin{rem}
For odd $n$, the vanishing of $w_{n-2}(\mathcal E)$ follows already from~(\ref{wh})
without any appeal to~(\ref{domrin}). Thus Lemma~\ref{norm} is true in that case
as well.
\end{rem}

\subsection{Dehn surgery}
\label{surgery}
Let $U\supset K$ be a tubular neighbourhood of a totally real embedded Klein bottle $K\subset\C^n$.
We consider two distinguished classes in the homology group $\H_{n-1}(\p U;\Z/2)$. Firstly, the fibre class $[\delta]$
generating the kernel of the inclusion homomorphism $\H_{n-1}(\p U;\Z/2)\to \H_{n-1}(\conj U;\Z/2)$
and, secondly, the class $[m^\sharp]$ of the $\C$-normal pushoff of $m$ introduced in \S\ref{ind}.

\begin{lem}[cf.\ \cite{N1}, Theorem 2.2]
\label{surg}
Consider a surgery $X=\conj U\cup_f (\C^n\setminus U)$ defined by a diffeomorphism
$f:\p U\to \p U$ such that
\begin{equation}
\label{dehn}
f_*[\delta]=[\delta]+[m^\sharp].
\end{equation}
If $n$ is even, then $K$ is homologically non-trivial in $X$. In particular, $\H_n(X;\Z/2)\ne 0$.
\end{lem}

\begin{proof}
As $n$ is even, we know that $\lk(K,m^\sharp)=1\mod 2$ by Lemma~\ref{viro} and the definition of the Viro index.
Since $\lk(K,\delta)=1\mod 2$ by definition, it follows that the sum $[\delta]+[m^\sharp]$
bounds a mod 2 chain in $\C^n-U$. By property~(\ref{dehn}), this chain and the $n$-ball bounded
by $\delta$ in $\conj U$ are glued into a mod 2 cycle in $X$ whose intersection index with $K$
is $1\mod 2$.
\end{proof}

\begin{lem}
\label{real}
If $X$ is orientable, then $\H_2(X;\R)=0$.
\end{lem}

\begin{proof}
Note first that $\H_2(X;\R)=\H_c^{2n-2}(X;\R)$ by Poincar\'e(--Lefschetz) duality.
Since $X\setminus K=\C^n\setminus K$, an inspection of the cohomology long exact sequences
$$
\dots\to \H^{2n-3}_c(K;\R)\to \H^{2n-2}_c(\C^n\setminus K;\R)\to
\H^{2n-2}_c(\C^n;\R)\cong 0
$$
$$
\dots\to \H^{2n-2}_c(X\setminus K;\R)\to
\H^{2n-2}_c(X;\R)\to\H^{2n-2}_c(K;\R)\cong 0
$$
shows that $\dim_\R \H^{2n-2}_c(X;\R)\le \dim_\R\H^{2n-3}_c(K;\R)$.
Thus, $\dim_\R\H_c^{2n-2}(X;\R)$ is zero for all $n\ge 3$ and does not exceed one for $n=2$.
In the latter case, however, it follows from Euler characteristic
additivity that the dimension of $\H^{2}_c(X;\R)$ is even and hence also equals zero.
\end{proof}

\subsection{Symplectic rigidity. Proof of the main result}
\label{mainproof}
If the surgery in Lemma~\ref{surg} were symplectic (i.\,e., there were a symplectic form on $X$
restricting to $\omega_0$ on $U$ and $\C^n\setminus\conj U$\,), then the conclusions
of Lemmas~\ref{surg} and~\ref{real} for an even $n$ would contradict the following result:

\begin{thm}[Eliashberg--Floer--McDuff~\cite{McD}, \cite{El}]
\label{rigid}
Let $(X,\omega)$ be a symplectic manifold symplectomorphic to $(\R^{2n},\omega_0)$, $n\ge 2$,
outside of a compact subset. Assume that $[\omega]$ vanishes on all spherical
elements in $\H_2(X;\R)$. Then $X$ is diffeomorphic to\/ $\R^{2n}$.
\end{thm}

\begin{rem}
If $n=2$, then $X$ is actually symplectomorphic to $(\R^4,\omega_0)$ by Gromov's
classical result~\cite{Gro}. Note, however, that we only need to know that $X$
must have the $\Z/2$-homology of the ball, which is proved in all dimensions
by a basic application of pseudoholomorphic curves (see~\cite{McD}, \S 3.8).
\end{rem}

Thus, to prove the main theorem it remains to show that for a {\it Lagrangian\/} embedding
of the Klein bottle $K$  there exists a symplectic surgery having property~(\ref{dehn}).
This can be done in all dimensions by the following elementary construction.

Represent $K$ as the quotient of $\R^n\setminus\{0\}$ by the $\Z$-action generated by the transformation
\begin{equation}
\label{act}
x\longmapsto 2\sigma(x),
\end{equation}
where $\sigma\in O_-(\R^n)$ is a reflection (in particular, $\sigma=\sigma^T=\sigma^{-1}$).
The cotangent bundle $T^*K$ is the quotient of $T^*(\R^n\setminus\{0\})\cong(\R^n_x\setminus\{0\})\times\R^n_y$
by the $\Z$-action generated by
\begin{equation}
\label{act1}
(x,y)\longmapsto (2\sigma(x),\frac12\sigma(y)).
\end{equation}
Note that the Riemannian metric $g=\frac{1}{\|x\|^2}{\sum dx_j^2}$ on $\R^n\setminus\{0\}$
is invariant with respect to~(\ref{act}) and equip $K$ with the induced metric. Note further
that the unit sphere bundle $ST^*(\R^n\setminus\{0\})\subset T^*(\R^n\setminus\{0\})$ with
respect to $g$ is the hypersurface $\{\|y\|^2=1/\|x\|^2\}$.

On $T^*(\R^n\setminus\{0\})$ with the zero section removed, consider the map
\begin{equation}
\label{map}
(x,y)\longmapsto (-y,x).
\end{equation}
Obviously, this map preserves the unit sphere bundle $ST^*(\R^n\setminus\{0\})$ and the canonical
symplectic form on $T^*(\R^n\setminus\{0\})$. Furthermore, it maps the orbits of the action~(\ref{act1}) into orbits.
Hence, it defines a symplectomorphism of $T^*K$ with the zero section removed that maps $ST^*K$ into itself.

Let us check that the action of the map (\ref{map}) on $\H_{n-1}(ST^*K;\Z/2)$ satisfies
condition~(\ref{dehn}). The fibre class $[\delta]$ is represented by the `vertical' $(n-1)$-sphere $\{x=\const, \|y\|=1\}$
and its image is obviously the class of the `horizontal' $(n-1)$-sphere $\{\|x\|=1,y=\const\}$.
Choose $m=\{\|x\|=1\}\subset K$ and $\nu_{m,K}(x)=x$. For any almost complex structure on $T^*K$ compatible
with the canonical symplectic form, the isotopy class of the $\C$-normal pushoff $m^\sharp=m+J\nu_{m,K}$
is the same as for the standard complex structure, i.\,e., it is given by the `diagonal' $(n-1)$-sphere $\{y=x, \|x\|=1\}\subset ST^*K$.
It follows immediately that the image of $[\delta]$ with respect to (\ref{map}) is $[\delta]+[m^\sharp]$, as required.

Finally, if $K$ is an embedded Lagrangian Klein bottle in a symplectic manifold,
we can identify its closed tubular neighbourhood $\conj U$ with the unit disc bundle $DT^*K$
by a conformally symplectic diffeomorphism and define the gluing map
$f:\p U\to\p U$ as the restriction of the symplectomorphism constructed above to $ST^*K$.\qed

\begin{rem}
\label{moresurg}
Replacing the action~(\ref{act}) by $x\mapsto 2x$, one obtains a completely analogous surgery
construction for the product $\Sph^{n-1}\times\Sph^1$. Further symplectic surgeries along
a Lagrangian Klein bottle or $\Sph^{n-1}\times\Sph^1$ can be defined by taking the gluing
map from the group generated by the map $f$ induced by~(\ref{map}) and the map $\tau$ induced
by the co-differential of the topologically non-trivial $g$-isometry $x\mapsto\frac{x}{\|x\|^2}$.
\end{rem}

\begin{rem}[Comparison with Luttinger surgery]
{\bf (i)}
In the case of the product $\Sph^{n-1}\times\Sph^1$, the surgeries found by Luttinger~\cite{Lu} for $n=2$
and by Borrelli~\cite{Bo} for $n=4$ and~$n=8$ correspond to the gluing maps $(f\circ\tau)^k$,
where $k\in\Z$ and the maps $f$ and $\tau$ are defined as in Remark~\ref{moresurg}.
{\bf (ii)}~%
The surgery used in~\cite{N1} in the case of the usual Klein bottle $K^2$ corresponds to
the gluing map $(f\circ\tau)^{-1}$. In the notation of~\cite{N1}, one has
$$
f(\phi,\psi,\theta)=(-\phi,\psi+\theta+\pi,-\theta-\pi)
\quad\mbox{ and }\quad
\tau(\phi,\psi,\theta)=(-\phi,\psi,-\theta-\pi)$$
so that $f\circ\tau=f_{0,-1}$.
{\bf (iii)}~%
There is an alternative description of these surgeries in terms of regluing Lefschetz pencils
via fibrewise symplectic Dehn twists (see the first draft of this paper, {\tt arxiv:0712.1760v1},
and the references therein).
\end{rem}

\begin{rem}[Totally real embeddings]
It is perhaps worth mentioning that totally real embeddings $K^n\hookrightarrow\C^n$
exist for all~$n$. Indeed, Lalonde~\cite{Lal} constructed Lagrangian immersions
$K^n\looparrowright\C^n$ that are regularly homotopic to embeddings. The existence
of totally real embeddings follows in this situation from Gromov's $h$-principle
(see, e.\,g.,~\cite{EM}, \S 19.3).
\end{rem}

\noindent
{\bf Acknowledgment.}
The author is grateful to the referee for useful comments.


{\small

\begin{thebibliography}{99}
\bibitem{Bo} V. Borrelli, {\it New examples of Lagrangian rigidity}, Israel J. Math. {\bf 125} (2001), 221--235.
\bibitem{Do} A. V. Domrin, {\it A description of characteristic classes of real submanifolds in complex manifolds
in terms of\/ $\R\C$-singularities}, Izv. Math. {\bf 59}:5 (1995), 899--918.
\bibitem{El} Y. Eliashberg, {\it On symplectic manifolds with some contact properties},
J. Differential Geom. {\bf 33} (1991), 233--238.
\bibitem{EP} Y. Eliashberg, L. Polterovich, {\it New applications of Luttinger's surgery},
Comment. Math. Helv. {\bf 69} (1994), 512--522.
\bibitem{EM} Y. Eliashberg, N. Mishachev, {\it Introduction to the $h$-principle},
Graduate Studies in Mathematics {\bf 48}, AMS, Providence, RI, 2002.
\bibitem{Gi} A. B. Givental, {\it Lagrangian imbeddings of surfaces and the open Whitney umbrella},
Functional Anal. Appl. {\bf 20}:3 (1986), 197--203.
\bibitem{Gro} M. Gromov, {\it Pseudoholomorphic curves in symplectic manifolds},
Invent. Math. {\bf 82} (1985), 307--347.
\bibitem{GM} L. Guillou, A. Marin (eds.), {\it \`A la recherche de la topologie
perdue}, Birkh\"auser, Boston, 1986.
\bibitem{La} H. F. Lai, {\it Characteristic classes of real manifolds
immersed in complex manifolds},
Trans. Amer. Math. Soc. {\bf 172} (1972), 1--33.
\bibitem{Lal} F. Lalonde, {\it Suppression lagrangienne de points doubles et rigidit\'e symplectique},
J. Differential Geom. {\bf 36} (1992), 747--764.
\bibitem{Lu} K. M. Luttinger, {\it  Lagrangian tori in $R\sp 4$},
J. Differential Geom. {\bf 42} (1995), 220--228.
\bibitem{McD} D. McDuff, {\it Symplectic manifolds with contact type boundaries},
Invent. Math. {\bf 103} (1991), 651--671.
\bibitem{McPh} R. MacPherson, {\it Generic vector bundle maps}, Dynamical Systems, Proc. Sympos.,
Univ. Bahia, Salvador, 1971, Academic Press, New York, 1973, pp. 165--175.
\bibitem{N} S. Nemirovski, {\it Lefschetz pencils, Morse functions,
and Lagrangian embeddings of the Klein bottle},
Izv. Math. {\bf 66}:1 (2002), 151--164.
\bibitem{N1} S. Nemirovski, {\it Homology class of a Lagrangian Klein bottle},
Preprint {\tt arxiv:math/0106122v4}, to appear in Izv. Math.
\bibitem{Sh} V. Shevchishin, {\it Lagrangian embeddings of the Klein bottle
and combinatorial properties of mapping class groups}, Preprint {\tt arxiv:0707.2085v1}.
\bibitem{We} S. M. Webster, {\it The Euler and Pontrjagin numbers of an $n$-manifold in $\C\sp n$},
Comment. Math. Helv. {\bf 60} (1985), 193--216.
\end{thebibliography}

}
\end{document}